\documentclass{amsart}
\usepackage{amsthm}
\usepackage{amstext}
\usepackage{amssymb}

\makeatletter
\numberwithin{equation}{section}
\numberwithin{figure}{section}
\theoremstyle{plain}
\newtheorem{thm}{Theorem}
  \theoremstyle{remark}
  \newtheorem*{rem*}{Remark}

\allowdisplaybreaks
\usepackage[sort&compress,square,numbers]{natbib}

\makeatother

\begin{document}

\title{A note on degenerate Bell numbers and polynomials}

\author{Tae Kyun Kim}
\address{Department of Mathematics, Kwangwoon University, Seoul 139-701, Republic of Korea}
\email{kimtk2015@gmail.com}

\author{Dae San Kim}
\address{Department of Mathematics, Sogang University, Seoul 121-742, Republic
of Korea}
\email{dskim@sogang.ac.kr}


\begin{abstract}
Recently, several authors have studied the degenerate Bernoulli and Euler
polynomials and given some intersting identities of those
polynomials. In this paper, we consider the degenerate Bell numbers
and polynomials and derive some new identities of those numbers
and polynomials associated with special numbers and polynomials. In addition,
we investigate some properties of the degenerate Bell polynomials
which are derived by using the notion of composita.
From our investigation, we
give some new relations between the degenerate Bell polynomials
and the special polynomials.
\end{abstract}

\maketitle
\global\long\def\relphantom#1{\mathrel{\phantom{{#1}}}}

\global\long\def\Zp{\mathbb{Z}_{p}}

\global\long\def\Bel{\mathrm{Bel}}

\global\long\def\acl#1#2{\left\langle \left.#1\right|#2\right\rangle }

\global\long\def\acr#1#2{\left\langle #1\left|#2\right.\right\rangle }

\section{Introduction}

As is well known, the ordinary Bernoulli polynomials are defined by
the generating function
\begin{equation}
\frac{t}{e^{t}-1}e^{xt}=\sum_{n=0}^{\infty}B_{n}\left(x\right)\frac{t^{n}}{n!},\quad\left(\text{see \cite{key-2,key-5}}\right).\label{eq:1}
\end{equation}
When $x=0$, $B_{n}=B_{n}\left(0\right)$ are called ordinary Bernoulli
numbers. From (\ref{eq:1}), we note that
\begin{equation}
B_{n}\left(x\right)=\sum_{l=0}^{n}\binom{n}{l}B_{l}x^{n-l},\quad\left(n\ge0\right),\quad\left(\text{see \cite{key-2,key-5}}\right).\label{eq:2}
\end{equation}

In \cite{key-3}, L. Carlitz considered the degenerate Bernoulli polynomials
which are given by the generating function
\begin{equation}
\frac{t}{\left(1+\lambda t\right)^{\frac{1}{\lambda}}-1}\left(1+\lambda t\right)^{\frac{x}{\lambda}}=\sum_{n=0}^{\infty}\beta_{n}\left(x\mid\lambda\right)\frac{t^{n}}{n!},\quad\left(\text{see \cite{key-1,key-6,key-8,key-9,key-12,key-17,key-18}}\right).\label{eq:3}
\end{equation}
When $x=0$, $\beta_{n}\left(\lambda\right)=\beta_{n}\left(0\mid\lambda\right)$
are called the degenerate Bernoulli numbers. These degenerate Bernoulli numbers
and polynomials are studied by several authors (see \cite{key-1,key-6,key-8,key-9,key-12,key-17,key-18}).

For $n\ge0$, the Stirling number of the first kind is defined as
\begin{align}
\relphantom{=} \left(x\right)_{n} & =x\left(x-1\right)\cdots\left(x-n+1\right)=\prod_{l=0}^{n-1}\left(x-l\right)\label{eq:4}\\
 & =\sum_{l=0}^{n}S_{1}\left(n,l\right)x^{l},\quad\left(\text{see \cite{key-10,key-16}}\right),\nonumber
\end{align}
and the Stirling number of the second kind is defined as
\begin{equation}
x^{n}=\sum_{l=0}^{n}S_{2}\left(n,l\right)\left(x\right)_{l},\quad\left(\text{see \cite{key-10,key-16}}\right).\label{eq:5}
\end{equation}

It is known that the generating functions of $S_{1}\left(n,l\right)$
and $S_{2}\left(n,l\right)$ are given by
\begin{equation}
\left(e^{t}-1\right)^{n}=n!\sum_{l=n}^{\infty}S_{2}\left(l,n\right)\frac{t^{l}}{l!},\label{eq:6}
\end{equation}
and
\begin{equation}
\left(\log\left(1+t\right)\right)^{n}=n!\sum_{l=n}^{\infty}S_{1}\left(l,n\right)\frac{t^{l}}{l!},\quad\left(\text{see \cite{key-10,key-16}}\right).\label{eq:7}
\end{equation}

The Bell polynomials (also called the exponential polynomial and denoted
by $\phi_{n}\left(x\right)$) are defined by the generating function
\begin{equation}
e^{x\left(e^{t}-1\right)}=\sum_{n=0}^{\infty}\Bel_{n}\left(x\right)\frac{t^{n}}{n!},\quad\left(\text{see \cite{key-4,key-7,key-13,key-14,key-15}}\right).\label{eq:8}
\end{equation}

It is not difficult to show that the first few of them are given by
\begin{align*}
\Bel_{0}\left(x\right) & =1,\quad\Bel_{1}\left(x\right)=x,\quad\Bel_{2}\left(x\right)=x^{2}+x,\quad\Bel_{3}\left(x\right)=x^{3}+3x^{2}+x,\\
\Bel_{4}\left(x\right) & =x^{4}+6x^{3}+7x^{2}+x,\quad\Bel_{5}\left(x\right)=x^{5}+10x^{4}+25x^{3}+15x^{2}+x,\\
\Bel_{6}\left(x\right) & =x^{6}+15x^{5}+65x^{4}+90x^{3}+35x^{2}+x,\quad\cdots.
\end{align*}
When $x=1$, $\Bel_{n}=\Bel_{n}\left(1\right)$ are called the Bell
numbers.

From (\ref{eq:8}), we can easily derive the following equation:
\begin{equation}
\Bel_{n}\left(x+y\right)=\sum_{l=0}^{n}\binom{n}{l}\Bel_{l}\left(x\right)\Bel_{n-l}\left(y\right),\quad\left(n\ge0\right),\label{eq:9}
\end{equation}
and
\begin{equation}
\frac{1}{e^{x}}\sum_{k=0}^{\infty}k^{n}\frac{x^{k}}{k!}=\sum_{k=0}^{n}S_{2}\left(n,k\right)x^{k}=\Bel_{n}\left(x\right),\quad\left(n\in\mathbb{N}\right).\label{eq:10}
\end{equation}

If we set $x=1$, then we obtain Dobi\'nski formula as follows:
\begin{equation}
\sum_{k=0}^{n}k^{n}\frac{1}{k!}=e\sum_{k=1}^{n}S_{2}\left(n,k\right)=e\Bel_{n},\quad\left(n\ge1\right),\label{eq:11}
\end{equation}
which is equivalent to
\begin{equation}
\Bel_{n}=\frac{1}{e}\sum_{k=0}^{n}k^{n}\frac{1}{k!},\quad\left(\text{see \cite{key-16}}\right).\label{eq:12}
\end{equation}

Let
\begin{align}
G\left(t,x\right) & =e^{x\left(e^{t}-1\right)}=\sum_{n=0}^{\infty}\Bel_{n}\left(x\right)\frac{t^{n}}{n!}\label{eq:13}
\end{align}

By differentiating $G\left(t,x\right)$ with respect to $t$, we get
\begin{align}
\relphantom{=} \sum_{n=0}^{\infty}\Bel_{n}\left(x\right)\frac{nt^{n-1}}{n!} & =\frac{d}{dt}G\left(t,x\right)=xe^{t}e^{x\left(e^{t}-1\right)}\label{eq:14}\\
 & =\sum_{n=0}^{\infty}\left(x\sum_{j=0}^{n}\Bel_{j}\left(x\right)\binom{n}{j}\right)\frac{t^{n}}{n!}.\nonumber
\end{align}

From (\ref{eq:14}), we have
\begin{equation}
\Bel_{n+1}\left(x\right)=x\sum_{j=0}^{n}\binom{n}{j}\Bel_{j}\left(x\right).\label{eq:15}
\end{equation}

In particular, for $x=1$, we get
\begin{equation}
\Bel_{n+1}=\sum_{j=0}^{n}\binom{n}{j}\Bel_{j}.\label{eq:16}
\end{equation}

Recently, several authors have studied the degenerate Bernoulli and Euler
polynomials and given some intersting identities of those
polynomials  (see \cite{key-1,key-6,key-8,key-9,key-12,key-17,key-18}).
In this paper, we consider the degenerate Bell numbers
and polynomials and derive some new identities of those numbers
and polynomials associated with special numbers and polynomials. In addition,
we investigate some properties of the degenerate Bell polynomials which are
derived by using the notion of composita. From our investigation, we
give some new relations between the degenerate Bell polynomials
and the special polynomials.

\section{Degenerate Bell polynomials and numbers}

Now, we consider the degenerate Bell polynomials which are given by
the generating funciton
\begin{equation}
\left(1+\lambda\right)^{\frac{x}{\lambda}\left(\left(1+\lambda t\right)^{\frac{1}{\lambda}}-1\right)}=\sum_{n=0}^{\infty}\Bel_{n,\lambda}\left(x\right)\frac{t^{n}}{n!}.\label{eq:17}
\end{equation}
When $x=1$, $\Bel_{n,\lambda}=\Bel_{n,\lambda}\left(1\right)$ are
called the degenerate Bell numbers.

From (\ref{eq:17}), we note that
\begin{align}
\sum_{n=0}^{\infty}\lim_{\lambda\rightarrow0}\Bel_{n,\lambda}\left(x\right)\frac{t^{n}}{n!} & =\lim_{\lambda\rightarrow0}\left(1+\lambda\right)^{\frac{x}{\lambda}\left(\left(1+\lambda t\right)^{\frac{1}{\lambda}}-1\right)}\label{eq:18}\\
 & =e^{\left(e^{t}-1\right)x}=\sum_{n=0}^{\infty}\Bel_{n}\left(x\right)\frac{t^{n}}{n!}.\nonumber
\end{align}

Thus, by (\ref{eq:18}), we get
\begin{equation}
\lim_{\lambda\rightarrow0}\Bel_{n,\lambda}\left(x\right)=\Bel_{n}\left(x\right),\quad\left(n\ge0\right).\label{eq:19}
\end{equation}

From (\ref{eq:17}), we can derive the following equation:
\begin{align}
 &\relphantom{=}  \sum_{n=0}^{\infty}\Bel_{n,\lambda}\left(x\right)\frac{t^{n}}{n!}\label{eq:20}\\
 & =\left(1+\lambda\right)^{\frac{x}{\lambda}\left(\left(1+\lambda t\right)^{\frac{1}{\lambda}}-1\right)}\nonumber \\
 & =\sum_{m=0}^{\infty}\left(\frac{\log\left(1+\lambda \right)}{\lambda}\right)^{m}\frac{1}{m!}x^{m}\left(\left(1+\lambda t\right)^{\frac{1}{\lambda}}-1\right)^{m}\nonumber \\
 & =\sum_{m=0}^{\infty}\left(\frac{\log\left(1+\lambda \right)}{\lambda}\right)^{m}\frac{x^{m}}{m!}\left(e^{\frac{1}{\lambda}\log\left(1+\lambda t\right)}-1\right)^{m}\nonumber \\
 & =\sum_{m=0}^{\infty}\left(\frac{\log\left(1+\lambda\right)}{\lambda}\right)^{m}x^{m}\sum_{k=m}^{\infty}S_{2}\left(k,m\right)\lambda^{-k}\frac{1}{k!}\left(\log\left(1+\lambda t\right)\right)^{k}\nonumber \\
 & =\sum_{k=0}^{\infty}\left(\sum_{m=0}^{k}\left(\frac{\log\left(1+\lambda\right)}{\lambda}\right)^{m}x^{m}S_{2}\left(k,m\right)\lambda^{-k}\sum_{n=k}^{\infty}S_{1}\left(n,k\right)\frac{\lambda^{n}}{n!}t^{n}\right)\nonumber \\
 & =\sum_{n=0}^{\infty}\left(\sum_{k=0}^{n}\sum_{m=0}^{k}\left(\frac{\log\left(1+\lambda\right)}{\lambda}\right)^{m}S_{2}\left(k,m\right)\lambda^{n-k}x^{m}S_{1}\left(n,k\right)\right)\frac{t^{n}}{n!}.\nonumber
\end{align}

By comparing the coefficients on the both sides of (\ref{eq:20}),
we obtain the following theorem.
\begin{thm}
\label{thm:1} For $n\ge0$, we have
\[
\Bel_{n,\lambda}\left(x\right)=\sum_{k=0}^{n}\sum_{m=0}^{k}\left(\frac{\log\left(1+\lambda\right)}{\lambda}\right)^{m}S_{1}\left(n,k\right)S_{2}\left(k,m\right)\lambda^{n-k}x^{m}.
\]

\end{thm}
Now, we observe that
\begin{align}
\left(1+\lambda\right)^{\frac{x}{\lambda}\left(1+\lambda t\right)^{\frac{1}{\lambda}}} & =\sum_{k=0}^{\infty}\left(\frac{x}{\lambda}\right)^{k}\frac{\left(\log\left(1+\lambda\right)\right)^{k}}{k!}\left(1+\lambda t\right)^{\frac{k}{\lambda}}\label{eq:21}\\
 & =\sum_{k=0}^{\infty}\left(\frac{x}{\lambda}\right)^{k}\left(\log\left(1+\lambda\right)\right)^{k}\frac{1}{k!}e^{\frac{k}{\lambda}\log\left(1+\lambda t\right)}\nonumber \\
 & =\sum_{k=0}^{\infty}\left(\frac{x}{\lambda}\right)^{k}\left(\log\left(1+\lambda\right)\right)^{k}\frac{1}{k!}\sum_{l=0}^{\infty}\left(\frac{k}{\lambda}\right)^{l}\frac{\left(\log\left(1+\lambda t\right)\right)^{l}}{l!}\nonumber \\
 & =\sum_{k=0}^{\infty}\left(\frac{x}{\lambda}\right)^{k}\frac{\left(\log\left(1+\lambda\right)\right)^{k}}{k!}\sum_{l=0}^{\infty}\left(\frac{k}{\lambda}\right)^{l}\sum_{n=l}^{\infty}S_{1}\left(n,l\right)\frac{\lambda^{n}}{n!}t^{n}\nonumber \\
 & =\sum_{k=0}^{\infty}\left(\frac{x}{\lambda}\right)^{k}\frac{\left(\log\left(1+\lambda\right)\right)^{k}}{k!}\sum_{n=0}^{\infty}\left(\sum_{l=0}^{n}k^{l}\lambda^{n-l}S_{1}\left(n,l\right)\right)\frac{t^{n}}{n!}\nonumber \\
 & =\sum_{n=0}^{\infty}\left\{ \sum_{k=0}^{\infty}\frac{x^{k}}{k!}\left(\frac{\log\left(1+\lambda\right)}{\lambda}\right)^{k}\sum_{l=0}^{n}k^{l}\lambda^{n-l}S_{1}\left(n,l\right)\right\} \frac{t^{n}}{n!}.\nonumber
\end{align}

Thus, by (\ref{eq:17}) and (\ref{eq:21}), we get
\begin{align}
 &\relphantom{=}  \left(1+\lambda\right)^{\frac{x}{\lambda}}\sum_{n=0}^{\infty}\Bel_{n,\lambda}\left(x\right)\frac{t^{n}}{n!}\label{eq:22}\\
 & =\sum_{n=0}^{\infty}\left\{ \sum_{k=0}^{\infty}\sum_{l=0}^{n}\frac{x^{k}}{k!}\left(\frac{\log\left(1+\lambda\right)}{\lambda}\right)^{k}k^{l}\lambda^{n-l}S_{1}\left(n,l\right)\right\} \frac{t^{n}}{n!}.\nonumber
\end{align}

Therefore, by Theorem (\ref{thm:1}) and (\ref{eq:22}), we obtain
the following theorem.
\begin{thm}
\label{thm:2} For $n\ge0$, we have
\begin{align*}
 & \left(1+\lambda\right)^{\frac{x}{\lambda}}\sum_{k=0}^{n}\sum_{m=0}^{k}\left(\frac{\log\left(1+\lambda\right)}{\lambda}\right)^{m}S_{1}\left(n,k\right)S_{2}\left(k,m\right)\lambda^{n-k}x^{m}\\
 & =\sum_{k=0}^{\infty}\sum_{l=0}^{n}\frac{x^{k}}{k!}\left(\frac{\log\left(1+\lambda\right)}{\lambda}\right)^{k}k^{l}\lambda^{n-l}S_{1}\left(n,l\right).
\end{align*}
\end{thm}
\begin{rem*}
$\:$
\begin{align*}
 &\relphantom{=}  e^{x}\sum_{m=0}^{n}S_{2}\left(n,m\right)x^{m}\\
 & =\lim_{\lambda\rightarrow0}\left(1+\lambda\right)^{\frac{x}{\lambda}}\sum_{k=0}^{n}\sum_{m=0}^{k}\left(\frac{\log\left(1+\lambda\right)}{\lambda}\right)^{m}S_{1}\left(n,k\right)S_{2}\left(k,m\right)\lambda^{n-k}x^{m}\\
 & =\lim_{\lambda\rightarrow0}\sum_{k=0}^{\infty}\sum_{l=0}^{n}\frac{x^{k}}{k!}\left(\frac{\log\left(1+\lambda\right)}{\lambda}\right)^{k}k^{l}\lambda^{n-l}S_{1}\left(n,l\right)\\
 & =\sum_{k=0}^{\infty}k^{n}\frac{x^{k}}{k!}.
\end{align*}

When $x=1$, we have
\begin{align}
 &\relphantom{=}  \left(1+\lambda\right)^{\frac{1}{\lambda}}\sum_{k=0}^{n}\sum_{m=0}^{k}\left(\frac{\log\left(1+\lambda\right)}{\lambda}\right)^{m}S_{1}\left(n,k\right)S_{2}\left(k,m\right)\lambda^{n-k}\label{eq:23}\\
 & =\sum_{k=0}^{\infty}\sum_{l=0}^{n}\frac{1}{k!}\left(\frac{\log\left(1+\lambda\right)}{\lambda}\right)^{k}k^{l}\lambda^{n-l}S_{1}\left(n,l\right).\nonumber
\end{align}

Note that
\begin{align*}
 & e\sum_{k=0}^{n}S_{2}\left(n,k\right)\\
 & =\lim_{\lambda\rightarrow0}\left(1+\lambda\right)^{\frac{1}{\lambda}}\sum_{k=0}^{n}\sum_{m=0}^{k}\left(\frac{\log\left(1+\lambda\right)}{\lambda}\right)^{m}S_{1}\left(n,k\right)S_{2}\left(k,m\right)\lambda^{n-k}\\
 & =\lim_{\lambda\rightarrow0}\sum_{k=0}^{\infty}\sum_{l=0}^{n}\frac{1}{k!}\left(\frac{\log\left(1+\lambda\right)}{\lambda}\right)^{k}k^{l}\lambda^{n-l}S_{1}\left(n,l\right)\\
 & =\sum_{k=0}^{\infty}k^{n}\frac{1}{k!},\quad\text{where }n\in\mathbb{N}.
\end{align*}

Now, we define the degenerate Stirling numbers of the second kind as
follows:
\begin{equation}
\left(e^{\frac{1}{\lambda}\log\left(1+\lambda t\right)}-1\right)^{n}=n!\sum_{m=n}^{\infty}S_{2}\left(m,n\mid\lambda\right)\frac{t^{m}}{m!}.\label{eq:23-1}
\end{equation}

Thus, from (\ref{eq:23-1}), we have

\begin{align}
 &\relphantom{=}  \left(1+\lambda\right)^{\frac{x}{\lambda}\left(\left(1+\lambda t\right)^{\frac{1}{\lambda}}-1\right)}=e^{\frac{x}{\lambda}\log\left(1+\lambda\right)\left(e^{\frac{1}{\lambda}\log\left(1+\lambda t\right)}-1\right)}\label{eq:24}\\
 & =\sum_{m=0}^{\infty}\frac{1}{m!}\left(\frac{x}{\lambda}\right)^{m}\left(\log\left(1+\lambda\right)\right)^{m}m!\sum_{k=m}^{\infty}S_{2}\left(k,m\mid\lambda\right)\frac{t^{k}}{k!}\nonumber \\
 & =\sum_{k=0}^{\infty}\left(\sum_{m=0}^{k}\left(\frac{\log\left(1+\lambda\right)}{\lambda}\right)^{m}x^{m}S_{2}\left(k,m\mid\lambda\right)\right)\frac{t^{k}}{k!}\nonumber \\
 & =\sum_{n=0}^{\infty}\left(\sum_{k=0}^{n}\left(\frac{\log\left(1+\lambda\right)}{\lambda}\right)^{k}x^{k}S_{2}\left(n,k\mid\lambda\right)\right)\frac{t^{n}}{n!}.\nonumber
\end{align}

By (\ref{eq:24}), we get
\begin{equation}
\Bel_{n,\lambda}\left(x\right)=\sum_{m=0}^{n}S_{2}\left(n,m\mid\lambda\right)\left(\frac{\log\left(1+\lambda\right)}{\lambda}\right)^{m}x^{m}.\label{eq:25}
\end{equation}

We observe that
\begin{align}
 &\relphantom{=}  \sum_{k=0}^{n}\left(\frac{\log\left(1+\lambda\right)}{\lambda}\right)^{k}S_{2}\left(n,k\mid\lambda\right)x^{k}\label{eq:26}\\
 & =\sum_{k=0}^{n}\sum_{m=0}^{k}\left(\frac{\log\left(1+\lambda\right)}{\lambda}\right)^{m}S_{1}\left(n,k\right)S_{2}\left(k,m\right)\lambda^{n-k}x^{m}\nonumber \\
 & =\sum_{m=0}^{n}\left(\sum_{k=m}^{n}S_{1}\left(n,k\right)S_{2}\left(k,m\right)\lambda^{n-k}\right)\left(\frac{\log\left(1+\lambda\right)}{\lambda}\right)^{m}x^{m}.\nonumber
\end{align}

Thus, by (\ref{eq:26}), we get
\begin{equation}
S_{2}\left(n,m\mid\lambda\right)=\sum_{k=m}^{n}S_{1}\left(n,k\right)S_{2}\left(k,m\right)\lambda^{n-k},\label{eq:27}
\end{equation}
where $0\le m\le n$.

Therefore, by (\ref{eq:26}) and (\ref{eq:27}), we obtain the following
theorem.\end{rem*}
\begin{thm}
\label{thm:3} For $n\in\mathbb{N}$, we have
\[
\Bel_{n,\lambda}\left(x\right)=\sum_{m=0}^{n}S_{2}\left(n,m\mid\lambda\right)\left(\frac{\log\left(1+\lambda\right)}{\lambda}\right)^{m}x^{m},
\]
where the degenerate Stirling numbers $S_{2}\left(n,m\mid\lambda\right)$ of the second kind have the expression
\[
S_{2}\left(n,m\mid\lambda\right)=\sum_{k=m}^{n}S_{1}\left(n,k\right)S_{2}\left(k,m\right)\lambda^{n-k},
\]
for $0\le m\le n$.
\end{thm}
Let us define the generating function of the degenerate Bell polynomials
as follows:
\begin{align}
G_{\lambda}\left(t,x\right) & =\left(1+\lambda\right)^{\frac{x}{\lambda}\left(\left(1+\lambda t\right)^{\frac{1}{\lambda}}-1\right)}\label{eq:28}\\
 & =\sum_{n=0}^{\infty}\Bel_{n,\lambda}\left(x\right)\frac{t^{n}}{n!}.\nonumber
\end{align}

Then, by (\ref{eq:28}), we get
\begin{align}
G_{\lambda}\left(t,x+y\right) & =\left(1+\lambda\right)^{\frac{x+y}{\lambda}\left(\left(1+\lambda t\right)^{\frac{1}{\lambda}}-1\right)}\label{eq:29}\\
 & =\left(1+\lambda\right)^{\frac{x}{\lambda}\left(\left(1+\lambda t\right)^{\frac{1}{\lambda}}-1\right)}\left(1+\lambda\right)^{\frac{y}{\lambda}\left(\left(1+\lambda t\right)^{\frac{1}{\lambda}}-1\right)}\nonumber \\
 & =\left(\sum_{m=0}^{\infty}\Bel_{m,\lambda}\left(x\right)\frac{t^{m}}{m!}\right)\left(\sum_{l=0}^{\infty}\Bel_{l,\lambda}\left(y\right)\frac{t^{l}}{l!}\right)\nonumber \\
 & =\sum_{n=0}^{\infty}\left(\sum_{m=0}^{n}\binom{n}{m}\Bel_{m,\lambda}\left(x\right)\Bel_{n-m,\lambda}\left(y\right)\right)\frac{t^{n}}{n!}.\nonumber
\end{align}

Therefore, by (\ref{eq:29}), we obtain the following theorem.
\begin{thm}
\label{thm:4} For $n\ge0$, we have
\[
\Bel_{n,\lambda}\left(x+y\right)=\sum_{m=0}^{n}\binom{n}{m}\Bel_{m,\lambda}\left(x\right)\Bel_{n-m,\lambda}\left(y\right).
\]

\end{thm}
By differentiating $G_{\lambda}\left(t,x\right)$ with respect to $t$, we get
\begin{align}
\sum_{n=1}^{\infty}\Bel_{n,\lambda}\left(x\right)\frac{nt^{n-1}}{n!} & =\frac{d}{dt}G_{\lambda}\left(t,x\right)\label{eq:30}\\
 & =\left(1+\lambda\right)^{\frac{x}{\lambda}\left(\left(1+\lambda t\right)^{\frac{1}{\lambda}}-1\right)}\frac{x}{\lambda}\log\left(1+\lambda\right)\left(1+\lambda t\right)^{-1+\frac{1}{\lambda}}\nonumber \\
 & =\left(\sum_{k=0}^{\infty}\Bel_{k,\lambda}\left(x\right)\frac{t^{k}}{k!}\right)\frac{x\log\left(1+\lambda\right)}{\lambda}\sum_{m=0}^{\infty}\left(\frac{1-\lambda}{\lambda}\right)_{m}\frac{\lambda^{m}t^{m}}{m!}\nonumber \\
 & =\sum_{n=0}^{\infty}\left(\sum_{k=0}^{n}\binom{n}{k}\Bel_{k,\lambda}\left(x\right)\left(1-\lambda\mid\lambda\right)_{n-k}\frac{x\log\left(1+\lambda\right)}{\lambda}\right)\frac{t^{n}}{n!},\nonumber
\end{align}
where $\left(x\mid\lambda\right)_{n}=x\left(x-\lambda\right)\cdots\left(x-\left(n-1\right)\lambda\right)$.

Therefore, by (\ref{eq:30}), we obtain the following theorem.
\begin{thm}
\label{thm:5} For $n\ge0$, we have
\[
\Bel_{n+1,\lambda}\left(x\right)=x\log\left(1+\lambda\right)^{\frac{1}{\lambda}}\sum_{k=0}^{n}\binom{n}{k}\Bel_{k,\lambda}\left(x\right)\left(1-\lambda\mid\lambda\right)_{n-k}
\]
where
\[
\left(x\mid\lambda\right)_{n}=x\left(x-\lambda\right)\cdots\left(x-\left(n-1\right)\lambda\right).
\]

\end{thm}
Note that
\begin{align}
 &\relphantom{=}  \sum_{k=0}^{\infty}\Bel_{k,\lambda}\left(x\right)\frac{t^{k}}{k!}\label{eq:31}\\
 & =\left(1+\lambda\right)^{\frac{x}{\lambda}\left(\left(1+\lambda t\right)^{\frac{1}{\lambda}}-1\right)}\nonumber \\
 & =\left(1+\lambda\right)^{-\frac{x}{\lambda}}\left(1+\lambda\right)^{\frac{x}{\lambda}\left(1+\lambda t\right)^{\frac{1}{\lambda}}}\nonumber \\
 & =\left(1+\lambda\right)^{-\frac{x}{\lambda}}e^{\left(1+\lambda t\right)^{\frac{1}{\lambda}}\log\left(1+\lambda\right)^{\frac{x}{\lambda}}}\nonumber \\
 & =\left(1+\lambda\right)^{-\frac{x}{\lambda}}\sum_{l=0}^{\infty}\frac{x^{l}}{l!}\left(\log\left(1+\lambda\right)^{\frac{1}{\lambda}}\right)^{l}\left(1+\lambda t\right)^{\frac{l}{\lambda}}\nonumber \\
 & =\left(1+\lambda\right)^{-\frac{x}{\lambda}}\sum_{l=0}^{\infty}\frac{x^{l}}{l!}\left(\log\left(1+\lambda\right)^{\frac{1}{\lambda}}\right)^{l}\sum_{k=0}^{\infty}\left(\frac{l}{\lambda}\right)_{k}\frac{\lambda^{k}t^{k}}{k!}\nonumber \\
 & =\left(1+\lambda\right)^{-\frac{x}{\lambda}}\sum_{k=0}^{\infty}\left(\sum_{l=0}^{\infty}\frac{x^{l}}{l!}\left(\frac{\log\left(1+\lambda\right)}{\lambda}\right)^{l}\left(l\mid\lambda\right)_{k}\right)\frac{t^{k}}{k!}.\nonumber
\end{align}

Therefore, by (\ref{eq:31}), we obtain the following theorem.
\begin{thm}
\label{thm:6} For $k\ge0$, we have
\[
\Bel_{k,\lambda}\left(x\right)=\left(1+\lambda\right)^{-\frac{x}{\lambda}}\sum_{l=0}^{\infty}\frac{x^{l}}{l!}\left(\frac{\log\left(1+\lambda\right)}{\lambda}\right)^{l}\left(l\mid\lambda\right)_{k}.
\]

\end{thm}
For $n\in\mathbb{N}$, we have

\begin{align}
\Bel_{n,\lambda}\left(x\right) & =\left(1+\lambda\right)^{-\frac{x}{\lambda}}\sum_{l=1}^{\infty}\frac{x^{l}}{l!}\left(\frac{\log\left(1+\lambda\right)}{\lambda}\right)^{l}\left(l\mid\lambda\right)_{n}\label{eq:32}\\
 & =\left(1+\lambda\right)^{-\frac{x}{\lambda}}x\sum_{l=0}^{\infty}\frac{x^{l}}{l!}\frac{\left(l+1\mid\lambda\right)_{n}}{l+1}\left(\frac{\log\left(1+\lambda\right)}{\lambda}\right)^{l+1}\nonumber \\
 & =\left(1+\lambda\right)^{-\frac{x}{\lambda}}\left(\frac{\log\left(1+\lambda\right)}{\lambda}\right)x\sum_{l=0}^{\infty}\frac{x^{l}}{l!}\frac{\left(l+1\mid\lambda\right)_{n}}{l+1}\left(\frac{\log\left(1+\lambda\right)}{\lambda}\right)^{l}\nonumber \\
 & =\left(1+\lambda\right)^{-\frac{x}{\lambda}}\left(\frac{\log\left(1+\lambda\right)}{\lambda}\right)x\sum_{l=0}^{\infty}\frac{x^{l}}{l!}\frac{1}{l+1}\nonumber\\
 &\relphantom{=}\times\sum_{k=0}^{n}\lambda^{n-k}S_{1}\left(n,k\right)\left(l+1\right)^{k}\left(\frac{\log\left(1+\lambda\right)}{\lambda}\right)^{l}\nonumber \\
 & =\left(1+\lambda\right)^{-\frac{x}{\lambda}}\left(\frac{\log\left(1+\lambda\right)}{\lambda}\right)x\sum_{k=0}^{n}\sum_{l=0}^{\infty}\frac{x^{l}}{l!}\lambda^{n-k}S_{1}\left(n,k\right)\nonumber\\
 &\relphantom{=}\times\left(l+1\right)^{k-1}\left(\frac{\log\left(1+\lambda\right)}{\lambda}\right)^{l}\nonumber \\
 & =\left(1+\lambda\right)^{-\frac{x}{\lambda}}\left(\frac{\log\left(1+\lambda\right)}{\lambda}\right)x\sum_{k=0}^{n}\lambda^{n-k}S_{1}\left(n,k\right)\nonumber\\
 &\relphantom{=}\times\sum_{l=0}^{\infty}\frac{x^{l}}{l!}\sum_{j=0}^{k-1}\binom{k-1}{j}l^{j}\left(\frac{\log\left(1+\lambda\right)}{\lambda}\right)^{l}\nonumber \\
 & =\left(1+\lambda\right)^{-\frac{x}{\lambda}}\left(\frac{\log\left(1+\lambda\right)}{\lambda}\right)x\sum_{k=0}^{n}S_{1}\left(n,k\right)\lambda^{n-k}\nonumber\\
 &\relphantom{=}\times\sum_{l=0}^{\infty}\frac{x^{l}}{l!}\sum_{j=0}^{k-1}\binom{k-1}{j}l^{j}\left(\frac{\log\left(1+\lambda\right)}{\lambda}\right)^{l}\nonumber \\
 & =\left(1+\lambda\right)^{-\frac{x}{\lambda}}e^{\frac{x}{\lambda}\log\left(1+\lambda\right)}\left(\frac{\log\left(1+\lambda\right)}{\lambda}\right)x\sum_{k=1}^{n}\sum_{j=1}^{k}S_{1}\left(n,k\right)\nonumber\\
 &\relphantom{=}\times\lambda^{n-k}\binom{k-1}{j-1}\Bel_{j-1}\left(\frac{\log\left(1+\lambda\right)}{\lambda}x\right).\nonumber
\end{align}

Therefore, by (\ref{eq:32}), we obtain the following theorem.
\begin{thm}
\label{thm:7} For $n\ge1$, we have
\[
\Bel_{n,\lambda}\left(x\right)=\frac{\log\left(1+\lambda\right)}{\lambda}x\sum_{k=1}^{n}\sum_{j=1}^{k}S_{1}\left(n,k\right)\lambda^{n-k}\binom{k-1}{j-1}\Bel_{j-1}\left(x\log\left(1+\lambda\right)^{\frac{1}{\lambda}}\right).
\]

\end{thm}
We observe that
\begin{align}
 &\relphantom{=}  \sum_{n=0}^{\infty}\frac{d}{dx}\Bel_{n,\lambda}\left(x\right)\frac{t^{n}}{n!}=\frac{d}{dx}\left(\left(1+\lambda\right)^{\frac{x}{\lambda}\left(\left(1+\lambda t\right)^{\frac{1}{\lambda}}-1\right)}\right)\label{eq:33}\\
 & =\frac{d}{dx}\left\{ e^{\frac{x}{\lambda}\left(\left(1+\lambda t\right)^{\frac{1}{\lambda}}-1\right)\log\left(1+\lambda\right)}\right\} \nonumber \\
 & =\frac{1}{\lambda}\left(\left(1+\lambda t\right)^{\frac{1}{\lambda}}-1\right)\log\left(1+\lambda\right)\left(1+\lambda\right)^{\frac{x}{\lambda}\left(\left(1+\lambda t\right)^{\frac{1}{\lambda}}-1\right)}\nonumber \\
 & =\left(1+\lambda t\right)^{\frac{1}{\lambda}}\left(\frac{\log\left(1+\lambda\right)}{\lambda}\right)\left(1+\lambda\right)^{\frac{x}{\lambda}\left(\left(1+\lambda t\right)^{\frac{1}{\lambda}}-1\right)}-\frac{\log\left(1+\lambda\right)}{\lambda}\left(1+\lambda\right)^{\frac{x}{\lambda}\left(\left(1+\lambda t\right)^{\frac{1}{\lambda}}-1\right)}\nonumber \\
 & =\left(\sum_{l=0}^{\infty}\left(1\mid\lambda\right)_{l}\frac{t^{l}}{l!}\right)\left(\sum_{m=0}^{\infty}\Bel_{m,\lambda}\left(x\right)\frac{t^{m}}{m!}\right)\left(\frac{\log\left(1+\lambda\right)}{\lambda}\right)-\frac{\log\left(1+\lambda\right)}{\lambda}\sum_{n=0}^{\infty}\Bel_{n,\lambda}\left(x\right)\nonumber \\
 & =\sum_{n=0}^{\infty}\left(\sum_{m=0}^{n}\binom{n}{m}\Bel_{m,\lambda}\left(x\right)\left(1\mid\lambda\right)_{n-m}\right)\frac{t^{n}}{n!}\frac{\log\left(1+\lambda\right)}{\lambda}-\frac{\log\left(1+\lambda\right)}{\lambda}\sum_{n=0}^{\infty}\Bel_{n,\lambda}\left(x\right)\frac{t^{n}}{n!}.\nonumber
\end{align}

By comparing the coefficients on the both sides of (\ref{eq:33}),
we obtain the following theorem.
\begin{thm}
\label{thm:8} For $n\ge1$, we have
\[
\frac{\lambda}{\log\left(1+\lambda\right)}\frac{d}{dx}\Bel_{n,\lambda}\left(x\right)=\sum_{m=0}^{n-1}\binom{n}{m}\Bel_{m,\lambda}\left(x\right)\left(1\mid\lambda\right)_{n-m}.
\]

\end{thm}

\section{Further remarks}

In \cite{key-11}, V. Kruchinin and D. Kruchinin introduced the notion
of composita in order to sutdy the coefficients of the powers of an
ordinary generating function and their properties.

Here we apply their technique to find an explicit expression of the
degenerate Bell polynomial $\Bel_{n,\lambda}\left(x\right)$. For
this, we first note that
\[
\sum_{n=0}^{\infty}\Bel_{n,\lambda}\left(x\right)\frac{t^{n}}{n!}=R\left(F\left(t\right)\right),
\]
where $R\left(t\right)=\left(1+\lambda\right)^{\frac{x}{\lambda}t}$,
$F\left(t\right)=\left(1+\lambda t\right)^{\frac{1}{\lambda}}-1$.

We recall from \cite{key-11} that the composita $G^{\Delta}\left(n,k\right)$of
the ordinary generating function $G\left(t\right)=\sum_{n=1}^{\infty}g\left(n\right)t^{n}$
is defined as the $n$th coefficient of $G\left(t\right)^{k}$. So
we have
\[
G\left(t\right)^{k}=\sum_{n=k}^{\infty}G^{\Delta}\left(n,k\right)t^{k}.
\]

Then it was noted in \cite{key-11} that
\[
G^{\Delta}\left(n,k\right)=\sum_{\lambda_{1}+\cdots+\lambda_{k}=n}g\left(\lambda_{1}\right)g\left(\lambda_{2}\right)\cdots g\left(\lambda_{k}\right),
\]
where the sum is over all compositions of the positive integer $n$
with $k$ parts.

In order to apply the following Theorem \ref{thm:9}, we need to determine
the composita $F^{\Delta}\left(n,k\right)$ of $F\left(t\right)$ and the coefficients $r\left(k\right)$ of $R\left(t\right)=\sum_{k=0}^{\infty}r\left(k\right)t^{k}$.
\begin{align*}
 & \sum_{k=0}^{\infty}r\left(k\right)t^{k}\\
 & =\left(1+\lambda\right)^{\frac{x}{\lambda}t}\\
 & =e^{\frac{x}{\lambda}\log\left(1+\lambda\right)t}\\
 & =\sum_{k=0}^{\infty}\frac{1}{k!}\left(\frac{x}{\lambda}\log\left(1+\lambda\right)\right)^{k}t^{k}\\
 & =\sum_{k=0}^{\infty}\frac{1}{k!}\left(\frac{\log\left(1+\lambda\right)}{\lambda}\right)^{k}x^{k}t^{k}.
\end{align*}

Then we obtain, for $k\ge0$,
\begin{equation}
r\left(k\right)=\frac{1}{k!}\left(\frac{\log\left(1+\lambda\right)}{\lambda}\right)^{k}x^{k}.\label{eq:34}
\end{equation}

Also, for $k\ge1$,
\begin{align*}
\sum_{n=k}^{\infty}F^{\Delta}\left(n,k\right)t^{n} & =F\left(t\right)^{k}\\
 & =\left(\left(1+\lambda t\right)^{\frac{1}{\lambda}}-1\right)^{k}\\
 & =\sum_{j=0}^{k}\binom{k}{j}\left(-1\right)^{k-j}\left(1+\lambda t\right)^{\frac{j}{\lambda}}\\
 & =\sum_{j=0}^{k}\binom{k}{j}\left(-1\right)^{k-j}\sum_{n=0}^{\infty}\left(j\mid\lambda\right)_{n}\frac{t^{n}}{n!}\\
 & =\sum_{n=0}^{\infty}\sum_{j=0}^{k}\left(-1\right)^{k-j}\binom{k}{j}\frac{\left(j\mid\lambda\right)_{n}}{n!}t^{n}.
\end{align*}

Hence, we have, for $k\ge1$,
\begin{align}
F^{\Delta}\left(n,k\right) & =\frac{1}{n!}\sum_{j=0}^{k}\left(-1\right)^{k-j}\binom{k}{j}\left(j\mid\lambda\right)_{n}\label{eq:35}\\
 & =\frac{1}{n!}\sum_{j=1}^{k}\left(-1\right)^{k-j}\binom{k}{j}\left(j\mid\lambda\right)_{n}.\nonumber
\end{align}

We now need the following theorem to get our result.
\begin{thm}[{\cite[Theorem 8]{key-11}}]
\label{thm:9} Suppose $F\left(t\right)=\sum_{n=1}^{\infty}f\left(n\right)t^{n}$,
with the composita $F^{\Delta}\left(n,k\right)$ and $R\left(t\right)=\sum_{n=0}^{\infty}r\left(n\right)t^{n}$.
Then, for the composition $A\left(t\right)=\sum_{n=0}^{\infty}a\left(n\right)t^{n}=R\left(F\left(t\right)\right)$,
the following holds:
\[
a\left(n\right)=\begin{cases}
r\left(0\right), & \text{for }n=0,\\
\sum_{k=1}^{n}F^{\Delta}\left(n,k\right)r\left(k\right), & \text{for }n\ge1.
\end{cases}
\]

\end{thm}
Now, the main result in this section follows from Theorem \ref{thm:9},
(\ref{eq:34}) and (\ref{eq:35}).
\begin{thm}
\label{thm:10} For all integers $n\ge0$, the degenerate Bell polynomial
$\Bel_{n,\lambda}\left(x\right)$ has the following expression
\[
\Bel_{n,\lambda}\left(x\right)=\begin{cases}
1, & \text{for }n=0,\\
\sum_{k=1}^{n}\left(\sum_{j=1}^{k}\frac{\left(-1\right)^{k-j}\binom{k}{j}}{n!k!}\left(\frac{\log\left(1+\lambda\right)}{\lambda}\right)^{k}\left(j\mid\lambda\right)_{n}\right)x^{k}, & \text{for }n\ge1.
\end{cases}
\]
\end{thm}
\bibliographystyle{amsplain}

\providecommand{\bysame}{\leavevmode\hbox to3em{\hrulefill}\thinspace}
\providecommand{\MR}{\relax\ifhmode\unskip\space\fi MR }
\providecommand{\MRhref}[2]{%
  \href{http://www.ams.org/mathscinet-getitem?mr=#1}{#2}
}
\providecommand{\href}[2]{#2}

\end{document}